\documentclass[a4paper,12pt]{article}
\usepackage{amssymb}
\usepackage{graphicx}
\usepackage{amsmath}
\usepackage{harvard}
\usepackage{lscape}

\input{tcilatex}

\begin{document}

\title{Incidence-based Estimates of Healthy Life Expectancy for the United
Kingdom: Coherence between Transition Probabilities and Aggregate Life Tables%
\thanks{%
Corresponding author:\ Ehsan Khoman, E-mail: e.khoman@niesr.ac.uk. This
paper is an adaptation of a more detailed report entitled, \textquotedblleft
Healthy Life Expectancy in the EU Member States\textquotedblright\ by the
same authors which has been recently published by the European Commission,
see www.enepri.org . We wish to thank A. Bebbington and J. Shapiro for
making available their estimates of the ordered probit model.}}
\author{Ehsan Khoman and Martin Weale \\
\\
National Institute of Economic and Social Research\\
2, Dean Trench Street,\\
London SW1P 3HE}
\maketitle

\begin{abstract}
Will the United Kingdom's ageing population be fit and independent, or
suffer from greater chronic ill health? Healthy life expectancy is commonly
used to assess this: it is an estimate of how many years are lived in good
health over the lifespan. This paper examines a means of generating
estimates of healthy and unhealthy life expectancy consistent with exogenous
population mortality data. The method takes population transition matrices
and adjusts these in a statistically coherent way so as to render them
consistent with aggregate life tables. It is applied to estimates of healthy
life expectancy for the United Kingdom.
\end{abstract}


\vspace*{0.5in} \textbf{\noindent Keywords}: Healthy Life Expectancy,
Least-squares Adjustment, Health State Transitions

\vspace*{0.5in}

\baselineskip0.27in

\newpage

\section{Introduction}

\qquad While it is plain that life expectancy has increased considerably
over the last thirty years or so in many advanced countries, it is much less
clear how healthy life expectancy has developed. Questions have therefore
arisen about the quality of life. Are we living longer but in worse health?
Are the increases in life expectancy at older ages because we are keeping
sick or disabled people alive longer or because we are saving people from
death but leaving them in states of disability? These are important
questions both for individuals and also for government policies on social
and health services provision for the elderly.

\qquad A shift in emphasis, from increasing survival to improving both the
length and quality of people's lives, has led to a greater policy interest
in the United Kingdom, and indeed in Europe as a whole in summary measures
of population health. The government projects that the overall number and
proportion of older people in the United Kingdom will rise significantly in
the coming decades (\citeasnoun{Wanless2002}). However, there is a debate
over whether these people will live longer, healthier lives, longer but more
disabled lives, or something in between. The UK Treasury's long-term
projections of the costs of an ageing population assume that the proportion
of life spent in long-term care will remain constant but acknowledge that
this is a cautious assumption and do not rule out an expansion of morbidity
for the United Kingdom (\citeasnoun{Treasury2004}).

\qquad A crucial question therefore is whether the proportion of life spent
in disability is expanding or decreasing. Existing data can be used to
support either case. While there have been clear rises in overall life
expectancy over time, there are concerns that not all years gained are in
good health and that the proportion of extra years lived are being spent in
ill-health (\citeasnoun{Bissett2002} and \citeasnoun{Breakwell2005}).
Therefore, the general consensus view in the academic community seems to be
that these trends reflect increased years of mild disability, and a decline
in severe disability (\citeasnoun{Bajekal2004} and \citeasnoun{Kelly2000}).

\qquad Existing calculations of healthy life expectancy are compiled from
the proportion of people reporting different health states (Sullivan's
method)- see appendix A for a comprehensive outline of Sullivan's method and
its uses. Health states of old people may reflect damage done in the past-
such as injuries sustained by soldiers and civilians during the Second World
War. They may therefore be a poor reflection of the risks of poor health
faced by young people. \citeasnoun{Bebbington1996} therefore argues that
healthy life expectancy should be calculated from the incidence of poor
health rather than its prevalence. In terms of acceptability, the fact that
transitions explicitly include return from poor to better states is
important. This enables a distinction to be made between a recurrent health
condition which allows for recovery, and one of steady decline to death.
Estimates of transition rates can be used for the prediction of lifetime
risk to individuals of particular states of ill-health, whereas prevalence
based measures cannot do this.

\qquad In order to produce measures of healthy life expectancy on this
basis, information is needed on transition matrices between different health
states. Such information may be available from household panel surveys such
as the European Community Household Panel (ECHP), carried out in the fifteen
countries of the pre-2005 European Union between 1994 and 2001. However,
such surveys are typically conducted on relatively small populations, and,
without further attention, the estimates of healthy and unhealthy life
expectancy generated by them are unlikely to be consistent with life tables
constructed from population mortality data.

\qquad In this paper we draw on a study of annual probabilities of
transition between health states by \citeasnoun{Bebbington2005}. We describe
a means of generating estimates of healthy and unhealthy life expectancy
consistent with exogenous population mortality data. The method takes
population transition matrices estimated from the ordered probit equations
in \citeasnoun{Bebbington2005} and adjusts these in a statistically coherent
way so as to render the transition matrices consistent with the mortality
data. It is applied to estimates of healthy life expectancy for the United\
Kingdom.

\qquad Since, health expectancy is a complex, multi-faceted concept, this
paper essentially aims to analyse the dynamics of health relating to the
transition of health states in the ECHP data. This paper answers the
following two questions. First, what the probability that an individual will
be in the same health state next year, be free of disability, be in worse
health or be dead? Secondly, what is the expected time spent in each health
state given that an individual is initially in a given health category?

\section{Data and Methodology}

\subsection{The ECHP}

\qquad This paper draws on the results of \citeasnoun{Bebbington2005}
presented in appendix B. They make use of the ECHP, the major innovative
attempt at a harmonised household (longitudinal) panel across the member
states of the European Union. The ECHP is essentially a standardised
multi-purpose annual longitudinal survey carried out between 1994 (wave 1)
to 2001 (wave 8) on each member state. Three characteristics make the ECHP a
unique source of information. These are (i) its multi-dimensional coverage
of a range of topics simultaneously; (ii) a standardised methodology and
procedures yielding comparable information across countries; and (iii) a
longitudinal or panel design in which information on the same set of
households and persons is gathered to study changes over time at the micro
level.

\qquad \citeasnoun{Bebbington2005} modelled the annual probabilities of
transition between health states for the EU\ member states including the
United Kingdom using pooled ordered probit equations from the ECHP. Separate
formulae were used for people above and below 65. Here we focus on the
results for the United Kingdom.

\subsection{Choice of Health Measures}

\qquad From the range of health status variables available in the ECHP, two
in particular were chosen. These are self-assessed health (SAH) (indicator
PiH001) and the existence of a chronic health or disability problem (PiH002)
combined with the degree of hampering health (HH) (PiH003).

\subsubsection{Self-Assessed Health}

\qquad In the ECHP users' database (UDB), self-assessed health (SAH) is
asked as `Please think back over the last 12 months about how your health
has been. Compared to people of your own age, would you say that your health
has on the whole been (i) excellent; (ii) good; (iii) fair; (iv) bad; or (v)
very bad? (PiH001)'. \citeasnoun{Bebbington2005} took the decision, after
considering the responses to PiH001 to combine `bad' and `very bad' health
states. Although this may remove some potential information, it avoids a
serious problem arising from the small numbers found in the worst category
in even the highest age groups.\ Therefore we can think of death as a fifth
state ranked below bad/very bad health.

\subsubsection{Hampering Health Condition}

\qquad The second measure of health is derived from the hampering health
(HH) condition. This indicator derives from two questions. Firstly, `Do you
have any chronic physical or mental health problem, illness or disability?
(PiH002)' and secondly, `Are you hampered in your daily activities by this
physical or mental health problem, illness or disability? (PiH003)'. The
three possible resulting states are (i) no such condition or a chronic
condition, but not hampered; (ii) hampered to some extent; or (iii) hampered
severely. Death is, as mentioned previously, an additional state.

\qquad \citeasnoun{Bebbington2005} ran into several serious problems
concerning the consistency and interpretation of the British data regarding
health, which are supplied to the ECHP as `clone' data from the British
Household Panel Survey (BHPS). A trial of three waves of parallel household
surveys, national and the ECHP, showed this was too much of a strain, with
high non-response rates, and as a result the sample size was reduced by
about a half from the fourth wave forwards. A conclusion from this is that
for HH, the category `to some extent' hampered was only used in the parallel
survey and then again in just wave 6 of the BHPS. The effect of this
seriously changed the distribution. In consequence %
\citeasnoun{Bebbington2005} made a decision to limit the analysis of the UK
sample by omitting the `to some extent' category, and on the evidence of the
UK parallel survey, results for this health definition will be incompatible
with other countries. As a result, only two health states were examined for
the HH measure of healthy life, namely, `no hampered condition' and
`hampered severely' with the additional absorbing state, death.

\qquad It is widely recognised that this indicator is less prone to
subjectivity than SAH and more immediately connected with disability,
dependency and a need for long-term care (\citeasnoun{Robine1998} and %
\citeasnoun{vandenBerg2001}). The European Commission considers this to be
an indicator for disability (\citeasnoun{Eurostat2001}). Also %
\citeasnoun{Bajekal2004} recently surveyed a variety of questions on
disability for the UK Department of Work and Pensions, and noted that a
similar census question which first made its appearance in 1991 had been
validated as a disability measure.

\section{Initial Transition Matrix Estimates}

\qquad For both of the domains distinguished, therefore, an ordered ranking
was generated running from the most healthiest state, `very good' for SAH,
and `not hampered to any degree' for HH, to the least favourable value, i.e.
death, the only absorbing state. The fact that the health states can be
ranked, a natural way to estimate transition probabilities as a function of
age and gender is by fitting an ordered probit model. Here we draw on the
models estimated by \citeasnoun{Bebbington2005} as mentioned above (see also %
\citeasnoun{Contoyannis2004a}), one for men and women aged under sixty-five
and the other for those aged sixty-five or older (Appendix B presents
estimates of the ordered probit equations for the United Kingdom derived
from \citeasnoun{Bebbington2005}).

\qquad The underlying probit function applied follows %
\citeasnoun{Wooldridge2002} and was used for example by %
\citeasnoun{Contoyannis2004b} in a similar analysis of health transitions
with the ECHP. A modelling approach to estimating transitions that makes use
of the latent variable specification can be written as%
\begin{equation}
h_{i}^{\ast }=\beta _{k}+x_{i}^{\prime }\cdot \gamma _{k}+e_{i}
\end{equation}%
where $\mathit{h}_{i}^{\ast }$ is some underlying continuous latent variable
for the$~\mathit{i}$th individual that underlies reported SAH and HH; $\beta
_{k}$ is a constant depending on the starting health state $k$; $e_{i}$
denotes a random, independently distributed component following a normal $N$%
(0,1) distribution. The variable $x_{i}$ is a vector of covariates including
age and gender coefficients and $\gamma _{k}$ a vector of parameters, which
again are assumed specific to the starting health state. If there is a
general trend, \citeasnoun{Bebbington2005} suggest that gender coefficients,
applying to women, tend to be positive at initial good states of health,
negative at bad states of health. This implies that women are more likely to
decline from good states of health, but men are more likely to decline or
die once in a bad state of health. \citeasnoun{Bebbington2005} have also
argued that it is plausible to drop the time-dependence $t$ in the present
case, and pool across waves 1 to 8 (i.e. 1994-2001), since there is no
discernible evidence of trends in the transitions. Since $\mathit{h}%
_{i}^{\ast }$ is not observed, \citeasnoun{Bebbington2005} in effect
partition it into the observed states, $\mathit{h}_{i}$,\ by a set of
unknown cut points, $\alpha $, (or threshold parameters), such that%
\begin{equation}
h_{i}=j\text{ if }\alpha _{j,k}<\mathit{h}_{i}^{\ast }\leq \alpha
_{j+1,k}~,~j=1,...,n
\end{equation}%
where $\alpha _{0}=-\infty $; $\alpha _{j}\leq \alpha _{j+1}$\ and $\alpha
_{J}=\infty $. Thus each observed health state corresponds to a value range
within the unobserved, latent distribution for health, such that the entire
range of the distribution is covered by one and only one health state. The
transition probabilities derive from the conditional distribution of $%
h_{i,t+1}$ given the state $k$ at time $t$:%
\begin{equation}
P(h_{i,t+1}=j\mid k)=\Phi (\alpha _{j+1,k}-\beta _{k})-\Phi (\alpha
_{j,k}-\beta _{k})
\end{equation}%
where $\Phi $ denotes the cumulative standardised normal distribution.

\qquad From these probit equations we calculate transition matrices as a
function of age and gender. We denote these \textbf{M}$_{\mathbf{0}}$ to 
\textbf{M}$_{\mathbf{99}}$\textbf{. }\ For an initial population vector 
\textbf{x}$_{i}$ whose $\emph{j}$th element, $x_{ij}$ shows the number of
people in health state $j$ on their $i$th birthday. It then follows that $%
\mathbf{x}_{i+1}=\mathbf{M}_{i}\mathbf{x}_{i}.$ If we denote by \textbf{i} a
vector of 1s with length equal to the number of health states, then from an
initial population \textbf{x}$_{0}$ the proportion surviving to their $i+1$%
st birthday is given as%
\begin{equation}
s_{i}=\frac{\mathbf{i}^{\prime }\Pi _{k=0}^{i-1}\mathbf{M}_{k}\mathbf{x}_{0}%
}{\mathbf{i}^{\prime }\mathbf{x}_{0}}
\end{equation}%
while we denote the proportion surviving to their $i$th birthday in the life
table as $s_{i}^{\ast }.$ Given our initial estimates of the transition
matrices, we wish to find new transition matrices, \textbf{M}$_{k}^{\ast }$
such that%
\begin{equation}
s_{i}^{\ast }=\frac{\mathbf{i}^{\prime }\Pi _{k=0}^{i-1}\mathbf{M}_{k}^{\ast
}\mathbf{x}_{0}}{\mathbf{i}^{\prime }\mathbf{x}_{0}}
\end{equation}%
where the $\mathbf{M}_{k}^{\ast }$ are reasonably close to the initial
estimates \textbf{M}$_{k}$. It is obvious that $s_{i}$ can be driven to $%
s_{i}^{\ast }$ only by adjusting the transition matrices \textbf{M}$_{k}$
where $k\leq i-1$. But an adjustment to one of these matrices has
implications for $s_{m}$ for all $m>i$. Thus, although it is obviously
possible to address the problem sequentially, it is unlikely that sequential
adjustment will offer the most satisfactory solution.

\section{A Least-Squares Approach}

\qquad Following from \citeasnoun{Deming1940} who first proposed the use of
a proportional fitting procedure to estimate cell probabilities in a
contingency table subject to certain marginal constraints, we set out here a
least-squares solution to the problem of adjusting the transition matrices
in order for them to be consistent with exogenous mortality data. Our
approach however differs somewhat from the original methods of %
\citeasnoun{Deming1940}, and that later extended by %
\citeasnoun{Friedlander1961}, in that while they concentrated on a linear
model in which the solution is derived exactly, the constraints that we face
are non-linear functions of the transition probabilities.

\qquad We denote by the vector \textbf{n}$_{k}$ the vector constructed from
the four columns of transition matrix \textbf{M}$_{k}$ stacked in order and
further consider the vector%
\begin{equation}
\mathbf{n}^{0}\mathbf{=}\left[ 
\begin{array}{c}
\mathbf{n}_{0} \\ 
... \\ 
\mathbf{n}_{k} \\ 
... \\ 
\mathbf{n}_{99}%
\end{array}%
\right] 
\end{equation}%
\qquad We write the vector of survival proportions generated by the vector $%
\mathbf{n}$ as $\mathbf{s(n)}$ with its $i$th element $s_{i}(\mathbf{n}%
)=s_{i}$ and the observed survival proportions as $\mathbf{s}^{\ast }.$ We
then aim to find $\mathbf{n}^{\ast }=\mathbf{n}^{0}\mathbf{+\Delta n}$ to
minimise 
\begin{equation}
\frac{1}{2}\mathbf{\Delta n}^{\prime }\mathbf{V}^{-1}\mathbf{\Delta n}+%
\mathbf{\lambda }\left\{ \mathbf{s}^{\ast }-\mathbf{s}\left( \mathbf{n}%
^{0}+\Delta \mathbf{n}\right) \right\} 
\end{equation}%
where $\mathbf{V}^{-1}$\ is a weighting matrix with $V_{ij}$ indicating the $%
i$th row and $j$th column of $\mathbf{V}$ with $n_{i}$ the $i$th element of $%
\mathbf{n}^{0}$. We set $V_{ii}=n_{i}^{2}$ and $V_{ij}=0$ ($i\neq j$) .\
Differentiating with respect to the elements of $\mathbf{n}$ 
\begin{equation}
\mathbf{V}^{-1}\mathbf{\Delta n}-\left( \frac{\partial \mathbf{s}}{\partial 
\mathbf{n}}\right) ^{\prime }\mathbf{\lambda =0}
\end{equation}%
where $\frac{\partial \mathbf{s}}{\partial \mathbf{n}}$ denotes a matrix
whose $i$th row and $j$th column consists of $\frac{\partial s_{i}}{\partial
n_{j}}.$ This gives 
\begin{equation}
\mathbf{\Delta n}=\mathbf{V}\left( \frac{\partial \mathbf{s}}{\partial 
\mathbf{n}}\right) ^{\prime }\mathbf{\lambda }
\end{equation}

We also note that by applying a Taylor series expansion we have%
\begin{equation}
\mathbf{s}^{\ast }-\mathbf{s}\left( \mathbf{n}^{0}+\Delta \mathbf{n}\right)
\cong \mathbf{s}^{\ast }-\mathbf{s}\left( \mathbf{n}^{0}\right) -\left( 
\frac{\partial \mathbf{s}}{\partial \mathbf{n}}%
{\vert}%
_{\mathbf{n}^{0}}\right) \mathbf{\Delta n}
\end{equation}%
\qquad

Given that%
\begin{equation}
\mathbf{s}^{\ast }-\mathbf{s(n}^{0})-\left( \frac{\partial \mathbf{s}}{%
\partial \mathbf{n}}%
{\vert}%
_{\mathbf{n}^{0}}\right) \mathbf{\Delta n\cong 0}
\end{equation}%
\qquad The exogenous survival rates will be approximately delivered if%
\begin{equation}
\mathbf{s}^{\ast }\mathbf{-s}\left( \mathbf{n}^{0}\right) \cong \left( \frac{%
\partial \mathbf{s}}{\partial \mathbf{n}}%
{\vert}%
_{\mathbf{n}^{0}}\right) \mathbf{\Delta n}
\end{equation}

We then set $\frac{\partial \mathbf{s}}{\partial \mathbf{n}}%
{\vert}%
_{\mathbf{n}^{0}}=\mathbf{S}_{0}$ and $\lambda _{0}\mathbf{=}\left\{ \mathbf{%
S}_{0}\mathbf{VS}_{0}^{\prime }\right\} ^{-1}\left( \mathbf{s}^{\ast }%
\mathbf{-s}\left( \mathbf{n}_{0}\right) \right) $. Therefore 
\begin{equation}
\mathbf{\Delta n}_{0}\mathbf{=VS}_{0}^{\prime }\left\{ \mathbf{S}_{0}\mathbf{%
VS}_{0}^{\prime }\right\} ^{-1}\left( \mathbf{s}^{\ast }\mathbf{-s}\left( 
\mathbf{n}_{0}\right) \right)
\end{equation}%
This finalises the first stage of the iteration process.

\qquad We now put $\mathbf{n}^{1}=\mathbf{n}^{0}+\Delta \mathbf{n}^{0}$ and
seek to find a vector $\Delta \mathbf{n}^{1}\mathbf{\ }$to minimise 
\begin{equation}
\frac{1}{2}\left( \mathbf{\Delta n}^{0}+\Delta \mathbf{n}^{1}\right)
^{\prime }\mathbf{V}^{-1}\left( \mathbf{\Delta n}^{0}+\Delta \mathbf{n}%
^{1}\right) +\mathbf{\lambda }\left\{ \mathbf{s}^{\ast }-\mathbf{s}\left( 
\mathbf{n}^{0}+\Delta \mathbf{n}^{0}+\Delta \mathbf{n}^{1}\right) \right\}
\end{equation}

Thus, with $\frac{\partial \mathbf{s}}{\partial \mathbf{n}}%
{\vert}%
_{\mathbf{n}^{1}}=\mathbf{S}_{1}$. We then have%
\begin{equation}
\mathbf{V}^{-1}\left( \mathbf{\Delta n}^{0}+\Delta \mathbf{n}^{1}\right) -%
\mathbf{S}_{1}^{\prime }\mathbf{\lambda =}0
\end{equation}%
and approximately 
\begin{equation}
\mathbf{s}^{\ast }\mathbf{-s}\left( \mathbf{n}^{1}\right) \cong \mathbf{S}%
_{1}\mathbf{\Delta n}^{1}
\end{equation}

This then yields 
\begin{equation}
\mathbf{S}_{1}\left( \mathbf{\Delta n}^{0}+\Delta \mathbf{n}^{1}\right) =%
\mathbf{S}_{1}\mathbf{V\mathbf{S}_{1}^{\prime }\lambda }
\end{equation}%
whence we have 
\begin{equation}
\left( \mathbf{\Delta n}^{0}+\Delta \mathbf{n}^{1}\right) =\mathbf{VS}%
_{1}^{\prime }\left\{ \mathbf{S}_{1}\mathbf{VS}_{1}^{\prime }\right\}
^{-1}\left\{ \mathbf{S}_{1}\mathbf{\Delta n}^{0}+\mathbf{s}^{\ast }\mathbf{-s%
}\left( \mathbf{n}^{1}\right) \right\}
\end{equation}

A further increment $\Delta \mathbf{n}^{2}$ is chosen to satisfy

\begin{equation}
\mathbf{V}^{-1}\left( \mathbf{\Delta n}^{0}+\Delta \mathbf{n}^{1}+\Delta 
\mathbf{n}^{2}\right) -\mathbf{S}_{2}^{\prime }\mathbf{\lambda =0}
\end{equation}%
and approximately 
\begin{equation}
\mathbf{s}^{\ast }\mathbf{-s}\left( \mathbf{n}^{2}\right) \cong \mathbf{S}%
_{2}\mathbf{\Delta n}^{2}
\end{equation}%
giving 
\begin{equation}
\left( \mathbf{\Delta n}^{0}+\Delta \mathbf{n}^{1}+\Delta \mathbf{n}%
^{2}\right) =\mathbf{VS}_{2}^{\prime }\left\{ \mathbf{S}_{2}\mathbf{VS}%
_{2}^{\prime }\right\} ^{-1}\left\{ \mathbf{S}_{2}\left( \mathbf{\Delta n}%
^{0}+\mathbf{\Delta n}^{1}\right) +\mathbf{s}^{\ast }\mathbf{-s}\left( 
\mathbf{n}^{2}\right) \right\}
\end{equation}

A recursive algorithm can be constructed 
\begin{equation}
\mathbf{\Delta n}^{j}=\mathbf{VS}_{j}^{\prime }\left\{ \mathbf{S}_{j}\mathbf{%
VS}_{j}^{\prime }\right\} ^{-1}\left\{ \mathbf{S}_{j}\sum_{i=0}^{j-1}\mathbf{%
\Delta n}^{i}+\mathbf{s}^{\ast }\mathbf{-s}\left( \mathbf{n}^{j}\right)
\right\} -\sum_{i=0}^{j-1}\mathbf{\Delta n}^{i}
\end{equation}%
with $\mathbf{n}^{j}\mathbf{=n}^{0}+\sum_{i=0}^{j-1}\mathbf{\Delta n}^{i}$
and subsequently, for any $j$ $\frac{\partial \mathbf{s}}{\partial \mathbf{n}%
}%
{\vert}%
_{\mathbf{n}^{j}}=\mathbf{S}_{j}$. Since the least-squares minimand is
evaluated afresh at each value of $\mathbf{n}^{j}$ an optimum is reached as $%
\mathbf{\Delta n}^{j}$ converges towards zero and the iterations can be
stopped when it is close to zero as defined by an appropriate tolerance
level. The adjusted vector $\mathbf{n}^{j}$ provides the transition matrices
at the $j$th iteration and when these are consistent with observed survival
rates, so too will be the healthy and unhealthy life expectancies derived
from them.

\section{Application to the United Kingdom}

\qquad Healthy life expectancy is given as the probability of being in
either a `very good' or `good' state given the condition of being in a `very
good' health state to begin with for SAH. For HH, healthy life expectancy is
simply given as the probability of being in a `none/slight' state
conditional on the probability of being in a `none/slight' state initially.
The table below provides estimates of healthy and unhealthy life expectancy
using both SAH and HH for men and women at age sixty-five in the United
Kingdom averaged between the period 1994 (wave 1) to 2001 (wave 8). The
unadjusted estimates are derived from the transition probabilities computed
with the ordered probit equations prior to the alignment having taken place.
The data published by the Office of National Statistics (ONS) are derived
from interim life tables based on three adjacent years provided by the
Government Actuary's Department. The life tables from \citeasnoun{GAD2005}
are available for three year windows in which the central year was chosen as
the average, for example, the year 1994 was computed from the window years
1993 to 1995.

\begin{table}[h]\centering%

\begin{center}
$%
\begin{tabular}{|cc|c|cc|c|}
\hline
\multicolumn{2}{|c|}{} & \multicolumn{3}{|c|}{\textbf{Unadjusted Estimates}}
& \textbf{ONS Estimates} \\ \hline
\multicolumn{2}{|c|}{} & Life & \multicolumn{2}{|c|}{Healthy life} & Life \\ 
\multicolumn{2}{|c|}{} & expectancy & \multicolumn{2}{|c|}{expectancy} & 
expectancy \\ \hline
\multicolumn{2}{|c|}{} &  &  & \% of &  \\ 
\multicolumn{2}{|c|}{Years} & Years & Years & lifetime in & Years \\ 
\multicolumn{2}{|c|}{} &  &  & healthy life &  \\ \hline
Men & SAH & 16.4 & 9.6 & 58.5 & 15.2 \\ 
Women & SAH & 14.9 & 9.2 & 61.7 & 18.5 \\ 
Men & HH & 17.3 & 12.0 & 69.3 & 15.2 \\ 
Women & HH & 15.6 & 11.0 & 70.5 & 18.5 \\ \hline\cline{6-6}
\end{tabular}%
$
\end{center}

\caption{Life expectancy and healthy life expectancy estimates using both
SAH and HH at age 65 for men and women between 1994 and 2001 in the United Kingdom calculated from the transition probabilities\label{key}}%
\end{table}%

\qquad Table 1 clearly sets out the problem. Whilst the transition
probabilities were pooled over the eight waves of the ECHP and thus life
expectancy using the unadjusted estimates is taken as an average over the
eight years (i.e. 1994-2001), life expectancy from the official data were
computed by taking the average from each wave for the sample year of the
ECHP. Therefore, life expectancy calculated from the transition
probabilities given by the unadjusted estimates suggests an apparent
discrepancy with the official data. The clear conclusion that can be
identified from table 1 is that the unadjusted estimates do not appear to
deliver the results of life expectancy provided by the official data,
casting doubt on the use of the associated estimates of healthy life
expectancy.

\qquad By using the alignment process derived by means of an least-squares
approach (as discussed in section 4), healthy life expectancy consistent
with official life expectancy data was calculated for a time-series of eight
years between 1994-2001, which are given in tables 2 and 3 below for men and
women at age sixty-five in the United Kingdom. The results depend, of
course, on the assumed mix of health states at age sixty-five. We have
generated this using the adjusted transition probabilities from birth. The
health state mix at age sixty-five is insensitive to that used at birth to
start the process. ONS estimates, based on Sullivan's method, for the same
period are also presented in order to give a comparison with our alignment
results. As one would expect from these tables, both measures of healthy
life expectancy tend to increase steadily with time for men and women. The
implementation of the adjustment process has meant, of course, that life
expectancy figures are identical to the ONS estimates since they are derived
from mortality tables provided by the Government Actuary's Department
projections which are taken for cohorts aged sixty-five between 1981 and
2054 (\citeasnoun{GAD2005}). 
\begin{landscape}
\begin{table}[h] \centering%

\begin{center}
$%
\begin{tabular}{|c|ccc|ccc|ccc|}
\hline
& \multicolumn{3}{|c|}{\textbf{Self-Assessed Health - SAH}} & 
\multicolumn{3}{|c|}{\textbf{Hampering Health - HH}} & \multicolumn{3}{|c|}{%
\textbf{ONS Estimates}} \\ \hline
& Life & \multicolumn{2}{c|}{Healthy life} & Life & \multicolumn{2}{c|}{
Healthy life} & Life & \multicolumn{2}{c|}{Healthy life} \\ 
& expectancy & \multicolumn{2}{c|}{expectancy} & expectancy & 
\multicolumn{2}{c|}{expectancy} & expectancy & \multicolumn{2}{c|}{expectancy
} \\ \hline
&  &  & \% of &  &  & \% of &  &  & \% of \\ 
Year & Years & Years & lifetime in & Years & Years & lifetime in & Years & 
Years & lifetime in \\ 
&  &  & healthy life &  &  & healthy life &  &  & healthy life \\ \hline
1994 & 14.5 & 8.9 & 61.4 & 14.5 & 10.4 & 71.7 & 14.5 & 11.0 & 75.9 \\ 
1995 & 14.7 & 8.9 & 60.5 & 14.7 & 10.5 & 71.4 & 14.7 & 11.3 & 76.9 \\ 
1996 & 14.8 & 9.0 & 60.8 & 14.8 & 10.6 & 71.6 & 14.8 & ... & ... \\ 
1997 & 15.0 & 9.1 & 60.7 & 15.0 & 10.7 & 71.3 & 15.0 & 11.7 & 78.0 \\ 
1998 & 15.2 & 9.2 & 60.5 & 15.2 & 10.8 & 71.1 & 15.2 & ... & ... \\ 
1999 & 15.4 & 9.3 & 60.4 & 15.5 & 11.0 & 71.0 & 15.4 & 11.5 & 74.7 \\ 
2000 & 15.7 & 9.5 & 60.5 & 15.7 & 11.1 & 70.7 & 15.7 & ... & ... \\ 
2001 & 15.9 & 9.6 & 60.4 & 15.9 & 11.2 & 70.4 & 15.9 & 11.6 & 73.0 \\ 
Increase from & 1.4 & 0.7 & -1.0 & 1.4 & 0.9 & -1.3 & 1.4 & 0.6 & -2.9 \\ 
1994 to 2001 &  &  &  &  &  &  &  &  &  \\ \hline
\end{tabular}%
$
\end{center}

\caption{Life expectancy and healthy life expectancy estimates at age 65 for
men between 1994 and 2001 in the United Kingdom\label{keya}}%
\end{table}%
\end{landscape}

\begin{landscape}
\begin{table}[h] \centering%

\begin{center}
$%
\begin{tabular}{|c|ccc|ccc|ccc|}
\hline
& \multicolumn{3}{|c|}{\textbf{Self-Assessed Health - SAH}} & 
\multicolumn{3}{|c|}{\textbf{Hampering Health - HH}} & \multicolumn{3}{|c|}{%
\textbf{ONS Estimates}} \\ \hline
& Life & \multicolumn{2}{c|}{Healthy life} & Life & \multicolumn{2}{c|}{
Healthy life} & Life & \multicolumn{2}{c|}{Healthy life} \\ 
& expectancy & \multicolumn{2}{c|}{expectancy} & expectancy & 
\multicolumn{2}{c|}{expectancy} & expectancy & \multicolumn{2}{c|}{expectancy
} \\ \hline
&  &  & \% of &  &  & \% of &  &  & \% of \\ 
Year & Years & Years & lifetime in & Years & Years & lifetime in & Years & 
Years & lifetime in \\ 
&  &  & healthy life &  &  & healthy life &  &  & healthy life \\ \hline
1994 & 18.1 & 11.3 & 62.4 & 18.1 & 12.9 & 71.3 & 18.1 & 12.9 & 71.3 \\ 
1995 & 18.2 & 11.3 & 62.1 & 18.2 & 12.9 & 70.9 & 18.2 & 13.0 & 71.4 \\ 
1996 & 18.3 & 11.4 & 62.3 & 18.3 & 13.0 & 71.0 & 18.3 & ... & ... \\ 
1997 & 18.4 & 11.4 & 62.0 & 18.4 & 13.0 & 70.7 & 18.4 & 13.1 & 71.2 \\ 
1998 & 18.5 & 11.4 & 61.6 & 18.5 & 13.1 & 70.8 & 18.5 & ... & ... \\ 
1999 & 18.6 & 11.5 & 61.8 & 18.6 & 13.2 & 71.0 & 18.6 & 13.1 & 70.4 \\ 
2000 & 18.8 & 11.6 & 61.7 & 18.8 & 13.4 & 71.3 & 18.8 & ... & ... \\ 
2001 & 19.0 & 11.7 & 61.6 & 19.0 & 13.4 & 70.5 & 19.0 & 13.2 & 69.5 \\ 
Increase from & 0.9 & 0.4 & -0.8 & 0.9 & 0.4 & -0.8 & 0.9 & 0.3 & -1.8 \\ 
1994 to 2001 &  &  &  &  &  &  &  &  &  \\ \hline
\end{tabular}%
$
\end{center}

\caption{Life expectancy and healthy life expectancy estimates at age 65 for
women between 1994 and 2001 in the United Kingdom\label{keyb}}%
\end{table}%
\end{landscape}

\qquad The last row of tables 2 and 3 demonstrates how life expectancy and
healthy life expectancy estimates have changed over the eight year time
span. It appears to be the case that HH estimates of healthy life expectancy
are markedly higher than those given by SAH estimates for both men and
women. This could be due to a number of reasons, for instance, the health
categories of the two healthy life measures could be interpreted differently
by different individuals and hence therefore more people stating a healthy
state of wellbeing for the HH estimate. In other words, since the HH
definition of healthy life expectancy is much more wider (i.e. less health
categories and thus more chance of being placed in a healthy state), so that
many conditions count as healthy, relative to the much narrower SAH
definition where many people are classed as unhealthy, healthy life
expectancy using HH may give the impression that the time spent in healthy
life will be relatively long, and will tend towards total life expectancy if
very few people are counted as unhealthy. As an aside, it is interesting to
note that though the HH estimates of healthy life expectancy are higher than
that given by the SAH measure, the percentage of time spent in healthy life
for HH tends to be decreasing at a higher rate than that of SAH for both men
and women.

\qquad In general, although life expectancy has risen for both men and women
using both measures of healthy life, the percentage of the lifetime spent in
ill-health tends to be increasing for both men and women. This suggests that
people are now living to ages in which they are more likely to experience
chronic diseases and disability, supporting the expansion of morbidity
hypothesis whereby as life expectancy increases, older people become more
vulnerable to chronic diseases and spend more time in ill-health and thus a
higher proportion of people with health problems survive to an advanced age (%
\citeasnoun{Gruenberg1977} and \citeasnoun{Olshansky1991}).

\qquad The last three columns present data published by the ONS of life
expectancy and healthy life expectancy between 1994 and 2001. The ONS
defines healthy life expectancy (HLE) from the age-specific prevalence
(proportions) of the population (rather than in incidence terms) in healthy
and unhealthy conditions and age-specific mortality information. Data for
1996, 1998 and 2000 were not published by the Statistical Office.

\qquad The method applied (discussed in \citeasnoun{Kelly2000}) uses the
General Household Survey (GHS) to provide estimates of healthy life
expectancy using the Sullivan method. The GHS asks a similar \ question to
that of the SAH measure used in the ECHP; `Over the last 12 mths would you
say your health has on the whole been good, fairly good or not good?'. From
tables 2 and 3 it appears that the ONS estimates of healthy life expectancy
are somewhat higher for both men and women relative to our measures, except
for HH for women, in which our estimates are slightly higher. One
explanation of this discrepancy between our healthy life estimates and the
ONS data could well be due to the fact that our estimates are based on
incidence rates (i.e. represent current health conditions and can help
predict future health care requirements) whilst ONS figures are more
prevalence based (i.e. dependent on past history). Prevalenced based
measures may underestimate (or overestimate) health expectancy, because the
prevalence of ill-health at a given age in the population reflect the past
probabilities of becoming ill at each younger age (\citeasnoun{Jagger2001}).
Moreover, while both our findings and the ONS figures tend to support an
expansion in morbidity, it is clear that our results suggest a slight
increase in healthy life expectancy relative to the ONS estimates based on
the Sullivan method.

\qquad In sum, one apparent conclusion from the analysis appears to be that
though there is some variation in our measures of healthy life expectancy
and that of ONS estimates, the alignment procedure significantly reduces the
dispersion of healthy life expectancy for both men and women. This could
suggest that the unadjusted results derived from the probit equations may
appear to give inaccurate estimates of healthy life expectancy whilst by
adjusting the transition matrices to render them statistically coherent with
exogenous population mortality data tends to have produce much more precise
estimates of healthy life expectancy.

\section{Discussion and Conclusions}

\qquad Since this paper outlines a longitudinal health study different from
that performed using cross-sectional data and Sullivan's method, this has
meant we have the novel advantage of being able to take account of
transitions into and out of various health states over time for the United
Kingdom. This multistate approach has the advantage over Sullivan's method
of providing health expectancy estimates based on current rather than
historical morbidity prevalence rates. The multistate life tables of the
transition probabilities and the expected time spent in each health state
also provides a clearer basis on which to predict service needs.

\qquad The results of this paper lend themselves to support the expansion of
morbidity hypothesis; where the additional gains in life expectancy are
spent in bad health while the number of years spent in good health remains
constant. Also, our results do indeed point to a slower worsening in healthy
life expectancy than the ONS estimates based on Sullivan's method.
Therefore, whilst our results and ONS estimates appear to suggest an
expansion in morbidity, our findings propose an improvement in healthy life
expectancy relative to the Sullivan method.

\qquad However, it has to be recalled that when using healthy life
expectancy measures, such as, SAH and HH, estimates can change over time
simply due to changes in individuals' subjective perceptions rather than a
true deterioration or improvement in the population's health. Hence, since
SAH and HH are subjective measures, meanings attached by respondents to the
categories may have changed over time due to medical advances. Also, both
health measures differ between different subgroups of the population.
Therefore account must be taken for individual's interpretation of the
different health states which may be affected by individuals age, gender and
socio-economic circumstances. The same issue of perception and
interpretation do not apply to total life expectancy, hence, the difference
between quality and quantity health measures.\qquad 

\bibliographystyle{agsm}
\bibliography{hle}

\appendix

\section{Appendix: Estimation of Healthy Life Expectancy}

\qquad The estimation of healthy life expectancy is based on the concept of
a closed population within a given period of time, in this case, using the
ECHP data between 1994 and 2001. Thus, this population does not account for
immigration or emigration. At the end of the period in question, the
population can be partitioned into those who die within the period and those
who are still alive. Of those still alive, the majority are expected to be
healthy, and some are expected to be unhealthy. Hence, a model can be built
that measures the health status of individuals who are alive at the same
time it accounts for those who die in the period in question. This section
reviews the techniques used to incorporate healthy life expectancy, namely,
prevalence-based life tables (Sullivan's method) based on the prevalence of
disability that is a stock that is dependent on past history) and
incidence-based life tables (multistate method which can adjust to represent
current health problems). Many researchers have indeed commented on the
differences between the Sullivan and multistate methods (%
\citeasnoun{Bebbington1991} and \citeasnoun{vandeWater1995}).

\qquad On the whole, experience has shown that Sullivan's method can,
generally, be recommended for its simplicity, relative accuracy, ease of
interpretation and suitability for long-term trends and comparisons between
populations and subgroups. Yet although most empirical research has used
Sullivan's method, its limitations are now well understood. In particular,
Sullivan's method is not suitable for detecting recent abrupt changes in
trends, nor for estimating incidence rates, prognosis, or life-time risk. It
is therefore better in principle to base future estimates of health care
needs on the current incidence of ill-health, rather than on current
prevalence. Incidence rates provide estimates of \ the current state of
health needs, and thus offer more accurate forecast of future health care
needs. Hence, the reason to apply incidence based measures here to predict
precise measures of healthy life expectancy.

\subsection{Sullivan's Method}

\qquad Sullivan's method (see \citeasnoun{Sullivan1966} and %
\citeasnoun{Sullivan1971}) requires only a population life table (which can
be constructed for a population using the observed mortality rates at each
age for a given time period) and prevalence data for the health states of
interest. Such prevalence rates can be obtained readily from cross-sectional
health or disability surveys carried out for a population at a point in
time. Surveys of this type are carried out regularly in the United Kingdom,
both at the national (\citeasnoun{Robine1991} and \citeasnoun{Matthews2006})
and regional level (\citeasnoun{Congdon2006}), and indeed across the EU
member states (\citeasnoun{Robine2003} and \citeasnoun{Robine2004}). Its
interest lies in its simplicity, the availability of its basic data and its
independence of the size and age structure of the population. The health
status of a population is inherently difficult to measure because it is
often defined differently among individuals, populations, cultures, and even
across time periods.

\qquad The objective of the Sullivan method is essentially to calculate the
expected life expectancy of groups of individuals currently at specified
ages if they lived the rest of their lives experiencing the age-specific
mortality rates observed for the population at a specific time. Thus the
technique essentially uses the age-specific mortality to calculate the
proportion of individuals alive at the beginning of an age interval that die
before reaching the next age group. Hence, this technique is a powerful tool
for estimating the remaining years of life that a group of individuals can
expect to live once they reach a certain age. The procedure for calculating
Sullivan's method is outlined below:

1. For each age/gender group obtain the life table schedules and the
expectation of life for the year of interest. Then calculate%
\begin{equation}
_{n}L_{x}=e_{x}l_{x}-e_{x+n}l_{x+n}  \tag{A1}
\end{equation}%
where $_{n}L_{x}$is the conventional life table measure of the average
number of person years lived in the age interval $x$ an $x+n$ (alternatively
this may be calculated from mortality rates).

2. Obtain the ill-health rate $_{n}d_{x}$ in each age-group observed in a
survey or census. If they are excluded, add the numbers in communal
establishments catering for the sick and disabled. Calculate the average
number of persons aged $x$ to $x+n$ living without ill-health in each
age/gender group as%
\begin{equation}
_{n}LWD_{x}=_{n}L_{x}(1-_{n}d_{x})  \tag{A2}
\end{equation}

3. Calculate life expectancy without ill-health as%
\begin{equation}
HLE_{x}=(\dsum ~_{n}LWD_{x})/l_{x}  \tag{A3}
\end{equation}%
where the summation is from age $x$ upwards. Hence equation (A3) presents
the proportion of years lived in a healthy state.

\qquad However, given the overall usefulness of the Sullivan method, it is
better in principle to base future estimates on health care needs on the
current incidence of ill-health, rather than on current prevalence.
Prevalence of chronic health conditions is affected only by past history in
that it is seen as a stock variable reflecting past flows, rather than
current health risks (\citeasnoun{Robine1999}). For example, past wars may
continue to affect current disablement rates, as may the past state of
health care, as conditions such as polio and thalidomide illustrate.
Therefore, if public health is changing, present prevalence may be a poor
guide to the future. This is one reason why it is inadvisable simply to
project current average age-specific expenditure rates to predict future
long term care needs. Incidence is a better guide to the current state of
health needs, and hence to predictions of future health. In this case
though, the Sullivan health expectancy remains a meaningful indicator of the
state of health at a population, rather than prediction at an individual,
level.

\qquad Consequently, although Sullivan's method fails to be a good predictor
of changes in the years an individual can expect to live in healthy years,
it does remain a meaningful indicator of the state of health of a population
at a starting point in time. Hence, it reflects the healthy years an
individual can expect to live only if current patterns of prevalences apply
during an entire lifetime.

\subsection{The Multistate Method}

\qquad Although empirical research has mainly used Sullivan's method of
calculating healthy expectancies, the approach used here applies the
multistate life table method for calculating healthy life expectancy.
Multistate life table methods for calculating health expectancies were first
proposed by \citeasnoun{Rogers1990} and \citeasnoun{Crimmins1996} to take
into account reversible transitions between one health state and another.
This approach is theoretically attractive since it allows one to calculate
health expectancies for population subgroups in a specific health state at a
given age, for example, those in a `very good' health state at age
sixty-five, whereas the Sullivan method gives only the average health
expectancy for the entire population at a given age. Hence the multistate
method is based on incidence rates that represent current health conditions.
The procedure therefore carried out in this study which is outlined below
generalises the multistate life table, which analyses the transition from a
given health state to another state or to the absorbing state, death.

\qquad The approach applied here therefore provides the critical link
between information on mortality and information on the spectrum of
non-fatal health experiences among the living. As an alternative to %
\citeasnoun{Bebbington2005}, where the results were divided between under
sixty-five year olds and people aged sixty-five or older, an attempt was
made to compute gender specific values for all age groups between 0 and 99
for each Member State.

\qquad The initial stage of out model consisted of calculating transition
probabilities by constructing normal distributions from the $\alpha $
coefficients derived from the probit equations in \citeasnoun{Bebbington2005}%
, for each health state and for each of the two measures of health
expectancy. We denote by $\mathbf{M}^{i}$ the transition matrix for an
individual aged $i$. Each element $\mathbf{M}_{j,k}^{i}$ shows the
probability that an individual in health state $k$ in year $i$ will be in
health state $j$ in year $i+1$. So the transition probabilities for each
Member States are therefore given by%
\begin{eqnarray}
\mathbf{N}_{j,k}^{1} &=&\mathbf{M}_{j,k}^{1}  \TCItag{A4} \\
\mathbf{N}_{j,k}^{i+1} &=&\mathbf{M}_{j,k}^{i+1}\cdot \mathbf{N}_{j,k}^{i} 
\TCItag{A5}
\end{eqnarray}%
where $\mathbf{N}_{j,k}^{i}~$is the probability that an individual is state $%
j$ conditional on him or her being in state $k$ at birth.

\qquad The next step consisted of simply computing the expected time in each
health state given that the individual was in a specific health category to
begin with, as a function of age and gender. It is apparent for all the
countries examined that as the age of the individual increases the expected
time spent in good health deteriorates and the time spent in bad health or
dying rises. It should also be noted that although the figures are presented
for ages 0 to 99, the oldest age reported for any country is 91, so beyond
this point figures may be of doubtful value. In order to calculate expected
time spent in each of the health states, denoted by $\mathbf{Z}_{j,k}^{i}$,
we have%
\begin{eqnarray}
\mathbf{Z}_{j,k}^{99} &=&\mathbf{M}_{j,k}^{99}  \TCItag{A6} \\
\mathbf{Z}_{j,k}^{99-i} &=&\mathbf{M}_{j,k}^{99-i}\cdot \mathbf{Z}%
_{j,k}^{100-i}+\mathbf{Z}_{j,k}^{100-i}  \TCItag{A7}
\end{eqnarray}%
\qquad \qquad \qquad \qquad \qquad \qquad \qquad

\qquad Equations (A6) and (A7) therefore provide the basis for determining
the expected number of years that an individual will spend state $j$
conditional on him or her being in state $k$ to begin with for each men and
women in the United Kingdom. In order to conclude this section it is
worthwhile recalling that while the Sullivan method of calculating healthy
life expectancy is based on prevalence rates, i.e. the prevalence of
disability that is a stock that is dependent on past history, the multistate
method applied here is based on incidence rates and thus can adjust to
represent current health conditions.

\section{Appendix: Ordered Probit Equations used to Construct Transition
Probabilities}

\qquad The tables below provide estimates of health transition rates
estimated from the ordered probit equations given in %
\citeasnoun{Bebbington2005} for the United Kingdom. Standard errors of
coefficients are shown in brackets. * denotes coefficients (age, gender) not
statistically significant (5\% level). The tables also exclude admissions to
a health-care institution.

(a) People under 65

\begin{table}[h] \centering%

\begin{center}
\begin{tabular}{|c|cccccc|}
\hline
Initial Health & $\alpha _{1}$ & $\alpha _{2}$ & $\alpha _{3}$ & $\alpha
_{4} $ & Age (years) & Gender \\ \hline
Very Good & 0.264 & 1.490 & 2.221 & 3.143 & -0.001* & 0.078 \\ 
& (0.045) & (0.046) & (0.055) & (0.138) & (0.001) & (0.027) \\ \hline
Good & -0.779 & 1.064 & 2.097 & 3.444 & 0.002 & 0.108 \\ 
& (0.032) & (0.033) & (0.037) & (0.116) & (0.001) & (0.019) \\ \hline
Fair & -1.093 & 0.311 & 1.733 & 3.141 & 0.013 & -0.002* \\ 
& (0.053) & (0.050) & (0.054) & (0.085) & (0.001) & (0.029) \\ \hline
Bad/Very Bad & -1.284 & -0.246 & 0.699 & 2.880 & 0.019 & -0.107 \\ 
& (0.106) & (0.100) & (0.101) & (0.121) & (0.002) & (0.053) \\ \hline
\end{tabular}
\end{center}

\caption{Ordered probit formulae coefficients of transition
probabilities for self-reported health (SAH) from the ECHP (all waves, pooled)\label{key1}}%
\end{table}%

\begin{table}[h] \centering%

\begin{center}
\begin{tabular}{|c|cccc|}
\hline
Initial Health & $\alpha _{1}$ & $\alpha _{2}$ & Age (years) & Gender \\ 
\hline
None/Slight & 2.381 & 3.622 & 0.015 & 0.113 \\ 
& (0.067) & (0.080) & (0.001) & (0.037) \\ \hline
Severe & 0.336 & 3.229 & 0.022 & -0.217 \\ 
& (0.142) & (0.187) & (0.003) & (0.067) \\ \hline
\end{tabular}
\end{center}

\caption{Ordered probit formulae coefficients of transition
probabilities for hampering health condition (HH) from the ECHP (all waves, pooled)\label{key2}}%
\end{table}%

\newpage

(b) People 65 and over

\begin{table}[h] \centering%

\begin{center}
\begin{tabular}{|c|cccccc|}
\hline
Initial Health & $\alpha _{1}$ & $\alpha _{2}$ & $\alpha _{3}$ & $\alpha
_{4} $ & Age (years) & Gender \\ \hline
Very Good & 1.955 & 3.110 & 3.687 & 3.924 & 0.026 & 0.007* \\ 
& (0.664) & (0.658) & (0.634) & (0.614) & (0.009) & (0.078) \\ \hline
Good & 0.515 & 2.302 & 3.220 & 3.644 & 0.023 & 0.079* \\ 
& (0.323) & (0.326) & (0.32) & (0.310) & (0.004) & (0.045) \\ \hline
Fair & -0.629 & 0.705 & 2.131 & 2.962 & 0.017 & -0.076* \\ 
& (0.319) & (0.318) & (0.316) & (0.308) & (0.004) & (0.054) \\ \hline
Bad/Very Bad & -1.250 & -0.285 & 0.738 & 2.244 & 0.017 & -0.285 \\ 
& (0.525) & (0.506) & (0.505) & (0.495) & (0.007) & (0.089) \\ \hline
\end{tabular}
\end{center}

\caption{Ordered probit formulae coefficients of transition
probabilities for self-reported health (SAH) from the ECHP (all waves, pooled)\label{key3}}%
\end{table}%

\begin{table}[h] \centering%

\begin{center}
\begin{tabular}{|c|cccc|}
\hline
Initial Health & $\alpha _{1}$ & $\alpha _{2}$ & Age (years) & Gender \\ 
\hline
None/Slight & 3.977 & 4.882 & 0.040 & 0.025* \\ 
& (0.393) & (0.386) & (0.005) & (0.063) \\ \hline
Severe & 0.612 & 2.795 & 0.020 & -0.210 \\ 
& (0.503) & (0.497) & (0.007) & (0.098) \\ \hline
\end{tabular}
\end{center}

\caption{Ordered probit formulae coefficients of transition
probabilities for hampering health condition (HH) from the ECHP (all waves, pooled)\label{key4}}%
\end{table}%

\end{document}